\documentclass[a4paper,10pt]{article}

\usepackage[paper=letterpaper,margin=1in]{geometry}

\usepackage{amsmath,amssymb,amsfonts,epsfig,setspace,bigstrut,longtable,array,breqn,url,color}

\usepackage{authblk}

\usepackage{biblatex}
\title{A Triumvirate of AI Driven Theoretical Discovery}

\author[*]{Yang-Hui He}

\affil[1]{London Institute for Mathematical Sciences, Royal Institution, London W1S 4BS, UK}
\affil[2]{Merton College, University of Oxford, OX14JD, UK}
\affil[*]{hey@maths.ox.ac.uk}

\begin{document}
\maketitle

\vspace{10mm}

\begin{abstract}
Recent years have seen the dramatic rise of the usage of AI algorithms in pure mathematics and fundamental sciences such as theoretical physics.
This is perhaps counter-intuitive since mathematical sciences require the rigorous definitions, derivations, and proofs, in contrast to the experimental sciences which rely on the modelling of data with error-bars. In this Perspective, we categorize the approaches to mathematical discovery as ``top-down'', ``bottom-up'' and ``meta-mathematics'', as inspired by historical examples.
We review some of the progress over the last few years, comparing and contrasting both the advances and the short-comings in each approach. 
We argue that while the theorist is in no way in danger of being replaced by AI in the near future, the hybrid of human expertise and AI algorithms will become an integral part of theoretical discovery.
\end{abstract}


\noindent \textbf{Key points:} 
\begin{itemize}
    \item Bottom-up mathematics is building statements from axioms and definitions.

    \item Meta-mathematics is to treat all statements and their derivations as a language.
    
    \item Top-down mathematics is theoretical research guided by intuition, experience and pure data.

    \item Special consideration should be given to pure, theoretical, noiseless data, which can lead to profound conjectures.

    \item AI has made significant advances in all these approaches of mathematical and theoretical research, but the human expert will not be replaced any time soon.

    \item Humans and AI will work in tandem for theoretical discovery.
\end{itemize}



\section{Introduction}
In a recent conversation, a friend who is an eye doctor recalled that when he was a trainee ophthalmologist, his mentor was a world-leading specialist who could recognize a diseased eye at a single glance from a wall of images. 
Yet, when pressed, his mentor often could not give more precise an answer than dismissing it as gut-instinct.
This, we concluded, is the fundamental reason why Ariticial Intelligence (AI) is taking over so many disciplines of human endeavour. 

For practical purposes, the vast majority of human activity falls under the adage of ``if it ain't broke, don't fix it''.
A medic's primary goal is to cure the patient; understanding the mechanisms of the disease is secondary.
ChatGTP can perform as well as an average undergraduate student in an exam not because a large language model (LLM) can have true comprehension but because for all but the very top students in the world the goal is to pass the test, not to understand the material.
Even (or perhaps especially) in the creative arts, a good definition of genius would be the mastery of style (via supervised learning of preceding samples) and the deviation therefrom in an inexplicable way.
Deep neural networks can mimic so much of human activities because they more often than not simply requires black-boxes: learning from trial-and-error, and performing to within a margin of error.

This can never be so for the scientist.
The raison d'\^etre of this tiny percentage of humanity is to understand and to question.
Such a compulsion for {explainability} is especially true for the {\it theorist}.
Here and throughout, we will use the term {\it theoretical science} to include both (1) pure mathematics and (2) the development and testing of theories and hypotheses using mathematics, exemplified by theoretical physics.
In other words, we will adhere to the British academic convention and consider the likes of theoretical physics or theoretical computer science as sub-disciplines of mathematics.

In some sense, the desire for {\it interpretability} and explicability is the very reason why the field of AI-assisted theoretical discovery has been a relative late-comer to fundamental science.
While as far back as the 1990s experimental physicists at CERN, in searching for new particles, were amongst the first scientists to use AI \cite{perret1990new}, it wasn't until 2017-8 that machine learning techniques emerged in theoretical physics \cite{He:2017aed,Carifio:2017bov,Krefl:2017yox,Ruehle:2017mzq,Liu:2017dzi} and pure mathematics \cite{He:2018jtw,Bull:2018uow}.

The past five years has seen tremendous progress in AI-driven theoretical investigations.
We give but a few of the representative examples out of the hundreds of papers that have come to define this exciting field.
In theoretical physics, these have included
particle phenomenology from string theory
\cite{Demirtas:2018akl,Mutter:2018sra,Cole:2018emh,Otsuka:2020nsk,Ruehle:2020jrk,Larfors:2020ugo,Bies:2020gvf,Perez-Martinez:2021zjj,Abel:2021ddu},
establishing dictionaries between field theory and deep-learning \cite{Hashimoto:2018ftp,Koch:2019fxy,Erbin:2021kqf}, theoretical cosmology
\cite{Liu:2017dzi,Jinno:2018dek,Rudelius:2018yqi}, quantum field theories
\cite{Chen:2020dxg,Chen:2020jjw,Kaspschak:2020zws,Harvey:2021oue,Gupta:2022vhe,Lal:2023dkj}, and uncovering fundamental symmetries \cite{Chen:2020jjw,Gal:2020dyc,Krippendorf:2020gny,udrescu2020ai,cornelio2021ai,liu2021machine,lemos2023rediscovering}. In parallel, in pure mathematics, these have included algebraic geometry
\cite{He:2017aed,He:2018jtw,Constantin:2018hvl,Klaewer:2018sfl,Altman:2018zlc,He:2019vsj,Grimm:2019bey,Erbin:2020tks,Berman:2021mcw,Berglund:2021ztg,coates2023machine}, algebraic structures and representation theory \cite{He:2019nzx,Bao:2020nbi,davies2021advancing,Dechant:2022ccf,He:2023wwt}, symbolic algebra and computation \cite{zaremba2014learning,England_2018,lample2019deep,E_2020,peifer2020learning,dabelow2024symbolic}, differential and metric geometry \cite{Ashmore:2019wzb,Jejjala:2020wcc,Ashmore:2021ohf,Larfors:2022nep}, number theory \cite{Alessandretti:2019jbs,He:2020eva,He:2020kzg,He:2020qlg,Abel:2022wnt}, graph theory and combinatorics \cite{He:2020fdg,Bao:2021auj,Bao:2021ofk,Berglund:2021ztg}, as well as
knot theory \cite{Gukov:2020qaj,Craven:2020bdz,davies2021advancing}. As to the question of AI's role in our civilization,
there is a mixture of optimism \cite{cohen2018machine,hudson2019,xu2021artificial,georgescu2022machines,thiyagalingam2022scientific,gukov2024rigor,fink,bengio2024machine} and anxiety \cite{satariano2023governments}.

The purpose of this Perspective is to attempt a balanced overview of the advances that have been made in the last few years \cite{He:2018jtw,lu2023survey,zhang2023ai,williamson2023deep,he2023machine,gukov2024rigor,alexwilkins}, whilst bearing in mind
the concerns and limitations \cite{greinerpetter2019machines,Dolotin_2023,Fajardo_Fontiveros_2023,kolpakov2023impossibility}.
We shall
caution against both euphoric over-enthusiasm and unwarranted fear, and convey that the future of the mathematical sciences will be an inevitable and hopefully harmonious union between the human and AI.

\subsection{What is the AI Mathematician?}
Regarding the issue of AI-driven theoretical discovery, the natural question arises as to ``What is the AI mathematician''?
First, we need to ask a more basic question: ``What is mathematics?'', or perhaps on a more pragmatic level, ``What is a mathematician?''.
One could look at it in three ways \cite{He:2021oav}, on which we will expound in detail, in light of how AI has been instrumental to each over recent years.
Not delving into the depth of the philosophy of mathematics and the philosophy of science, one can loosely categorize as follows:

\begin{description}
    \item[I. Bottom-Up: ]
    One can think of mathematics as being built from foundational axioms, where all theorems and equations are constructed from the roots up using logic.
    This is what is known as Hilbert's Formalism Programme \cite{zach2007hilbert}.
    We will refer to this approach as ``bottom-up'' to reflect the rigorous nature of theoretical research.

    \item[II. Meta-Mathematics: ]
    Closely related to the bottom-up approach is the Logicism of Russell-Whitehead \cite{russell1910principia} and ultimately Wittgenstein \cite{wittgenstein}.
    We will think of this as ``Mathematics as Language'', where one considers any proposition as a set of symbols, led to by sequences of symbols that one calls {\it proof} or {\it derivation}.
    We will refer to this as meta-mathematics in the sense of looking at the problem from a distance, perhaps more as a linguist or computer-scientist.

    \item[III. Top-Down: ]
    The practicing theorist often {\it experiments} and {\it conjectures} before tackling a proof or derivation.
    This is somewhat in the spirit of
    Brouwer's intuitionist approach to mathematics \cite{brouwer1927intuitionistic} where one factors in the human element.
    We will refer to this as ``top-down'', where an over-arching view, based on experience and speculation, guides one towards a problem.

\end{description}
The above is the ``Triumvirate'' in the title.
We shall discuss how each has witnessed dramatic advancement in the last few years.

\section{Bottom-Up Mathematics}
Hilbert's 
``Wir mussen wissen. Wir werden wissen.'' (We must know. We will know.) is a famous declaration that should be considered in conjunction with the {\it Principia Mathematica} of Russell-Whitehead, in a tradition that dates back at least to debates with Frege \cite{frege1893grundgesetze}, or to the binary machines of Leibniz, or indeed to the {\it Elements} of Euclid.
The Programme of building up mathematical truths from the ground up received a fatal blow from the Undecidability and Incompleteness Theorems of G\"odel \cite{godel1931formal}, Church \cite{church1936unsolvable} and Turing \cite{turing1938computable} by the 1930s.

However, to quote Prof.~Minhyong Kim in a private communication, ``the practicing mathematician rarely contemplates whether your daily proposition is provable or not''. In other words, the space of decidable and interesting statements are so vast that one could first focus on these.
Thus, despite the logical impossibility of building all mathematical statements bottom-up, theorists certainly never stopped pursuing proofs for countless propositions.
This lead to the modern day answer to Russell-Whitehead's {\it Principia Mathematica} \cite{russell1910principia}: the Automated Theorem Proving  programme (ATP).
Arguably the first AI system for mathematics - or indeed, the first AI system - was the Logic Theory Machine \cite{newell1956logic} of Newell-Simon-Shaw of 1956, which was an early computer system \footnote{Interestingly, this was around the same time as the emergence of the first trainable neural network \cite{rosenblatt1958perceptron}.
} designed for, and succeeded in, proving a number of propositions of the {\it Principia}.

It has been some 70 years since this first AI-for-mathematics system and much progress has been made.
Over the second half of the twentieth century, it became clear that an increasing number of proofs of fundamental results in mathematics are impossible without the computer.
These have ranged from situations where key steps reduce to extensive brute-force computation, such as in the four-colour theorem, to more extreme circumstances where it takes longer than a human life-span to go through all the details, such as the classification of simple finite groups.
Dependence of the human theorist on machines have prompted such influential figures as Terrance Tao \cite{roberts2023ai} and International Congress of Mathematicians addresses \cite{icm} to seriously consider the future of mathematics.

While the first {\it proof-assistant} appeared in the 1970s,
\cite{de1994mathematical}, 
Isabelle/HOL \cite{nipkow2002isabelle}, Coq \cite{bertot2013interactive}, Agda \cite{wadler2018programming}, and Lean \cite{de2015lean} softwares are spear-heading the ATP programme in this century.
One notable direction well under way is the Xena project (\url{https://xenaproject.wordpress.com/}) headed by Prof.~Kevin Buzzard to formalize all (every statement and every step of proof) of undergraduate-level mathematics into Lean.
More recently and more non-trivially, Gowers, Green, Manners and Tao 
\cite{gowers2023conjecture}
used Lean's MathLib library to prove the Polynomial Freiman-Ruzsa conjecture.
Over several private conversations with Buzzard and Davenport, we are still far from having established anything close to a full database of all of contemporary mathematics in Coq or Lean format\footnote{Only earlier this year was a new project launched to formalize all the requisite pieces to Wiles' proof of Fermat's last theorem (UK EPSRC Grant EP/Y022904/1).}, let alone have AI automation on selecting correct proof strategies given a proposition or conjecture.
Nevertheless, given such a database \cite{openAI,metaAI}, there inevitably will be a plethora of research devoted to mining this data for new theories.
This brings us to the next point.

\section{Meta-Mathematics}
From the {\it Principia} to the advancement of computer science - or indeed, from Euclid's {\it Elements} or Galileo's {\it Il Saggiatore} - there has been a tradition of viewing mathematics as a language
\cite{ganesalingam2013language}.
Indeed, Natural Language Processing (NLP) is rooted in Turing's original proposal of his famous eponymous test \cite{turingcomputing}.
AI and the internet have propelled NLP to the era of the Large Language Model (LLM). Indeed, openAI's ChatGTP (or its counterparts such as google's Gemini) has passed the Turing Test \cite{biever2023chatgpt}.
The important point here, of course, is that LLM has no {\it understanding} of the underlying material, it is merely grouping together words in the right order based on large corpora of statistical samples.
The philosophy \footnote{Perhaps for this very reason of whether there is understanding of the underlying mathematics that we have chosen to call this direction ``meta-mathematics''.} of ``understanding'' aside, it is indubitable that ChatGTP has been transformative in mimicking human communication. 

One of the earliest experiments \cite{he2018hepth} on LLM for theory was the application of the Word2Vec \cite{mikolov2013efficient} neural network (perhaps the most basic LLM technique) applied to the titles of several sections of the ArXiv (\url{www.arxiv.org}), the most comprehensive repository for contemporary research in mathematics and theoretical physics.
Perhaps more interesting than in retrieving seemingly sensible linguistic identities (e.g., `string-theory + Calabi-Yau = M-theory + G2'), was a comparison with viXra (\url{www.viXra.org}), the repository of fringe ideas not accepted by main-stream science.
From the titles alone,  
one could significantly distinguish (from the confusion matrix)
different sub-fields of theoretical physics (high energy theory, high energy phenomenology, general relativity/quantum cosmology, etc) in arXiv, whereas in viXra this is not so. In other words, the syntax of proper theoretical science is more self-consistent than that of fringe science even at the level of titles\footnote{
If one enriched the data, and included abstracts, application of Word2Vec on papers in material science have uncovered new chemical reactions
\cite{tshitoyan2019unsupervised}.
}.

Today, Word2Vec has been superseded by transformers which are the preferred architectures for NLP, and the programme of {\it LLM for mathematics} is blooming
\cite{kim2021distilling,hutson2022ai,wang2023generative,romera2024mathematical,trinh2024solving,ahn2024large,frieder2024large}.
Notably, in parallel to the aforementioned OpenAI and MetaAI projects, DeepMind's AlphaGeo \cite{trinh2024solving} has recently been able to generate correct, human-understandable proofs for Olympiad level problems in Euclidean geometry.
It is clear, when and if we do have a linguistic database of all contemporary mathematics - which is certainly many decades in the future - the LLM approach of AlphaGeo on this vast data should produce new mathematics.

\section{Top-Down Mathematics}
Everything we have said so far has to do with building correct mathematical statements.
But frequently one has no idea {\it what} statement one should try to show.
Indeed, how does the practicing mathematician actually work?
In many ways, our papers are written backwards.
From day to day, we doodle on paper and on board, experimenting with ideas, mistakes, and expressions, until something sensible comes out.
Then, we go back and formalize with definitions followed by theorems and derivations that lead to logical conclusions.
Thus, journal papers in mathematics and theoretical physics look ``bottom-up'' even though the discovery process is quite the opposite.
Historically, this is even more true for some of the greatest theoretical discoveries.
Newton and Euler were freely manipulating formal expressions in calculus, centuries before a proper notion of analysis and convergence; Galois showed the unsolvability of the quintic by radicals by seeing the structures of permutation groups, before the definitions of groups and fields that we are taught today; theoretical physicists freely manipulate Feynman integrals to give results that agree with experiment to astounding accuracy, when we still do not have a mathematically rigorous formulation of quantum field theory
etc.

In a recent AI safety conference, a policy maker jokingly said to me that mathematicians are high on the list of jobs being replaced\footnote{The eye doctor mentioned earlier confessed to me that he will not take on any further trainees in diagnosis; that will soon be entirely taken over by AI.}
In her mind, a human mathematician is a bottom-up Logical Theory Machine, building sentences (proofs and derivations) from definitions.
In reality,
actual mathematical research is based on a combination of inspiration, intuition, and experience. To contrast with the almost dry narrative of ``bottom-up'', this almost fuzzy approach we will call ``top-down''.
Of course, at the end of the day, all statements must be rigorous and any fuzziness and inaccuracies must be distilled out (see nice recent Perspective \cite{gukov2024rigor}).
It is this direction which we now discuss in detail, with illustrative examples.

Perhaps contrary to common conception, an indispensable component to even the purest mathematical discovery is {\it data}.
This is not experimental data in the sense of, e.g., particle trajectories from CERN, with errors and variance, but results of classifications and computations, e.g., tables of characters of finite groups. This ``pure data'' is exact and without statistical error and shed light on the underlying theory.
To quote Vladimir Arnol'd, ``mathematics is the part of physics where experiments are cheap.'' Indeed, while AI has long been instrumental in the scientific big data revolution \cite{graham2013machine,hey2020machine}, its application to pure mathematical data is new \cite{He:2018jtw,He:2021oav}.

The best neural network of the nineteenth century is undoubtedly the brain of Gauss.
Confronted with the ancient problem  of finding patterns in primes which date back at least to Euclid, the sixteen-year-old defined the prime counting function which gives the number of prime numbers not exceeding a positive real $x \in \mathbb{R_+}$,
\begin{equation}
    \pi(x) = \#\{ p \leq x : p \mbox{ prime } \} \ .
\end{equation}
Gauss consulted tables available at the time and computed tens of thousands more (by hand!) and simply plotted $\pi(x)$.
With divine inspiration, he conjectured that $\pi(x) \sim x/\ln(x)$.
This profound observation had to wait for the establishment of complex analysis by Cauchy and Riemann in order to be proven by Hadamard and de la Vall\'ee Poussin at least 50 years later. It is now known as the Prime Number Theorem, a cornerstone of mathematics.

In the twentieth century, Birch and Swinnerton-Dyer plotted, using the earliest computers of the 1960s, ranks and other quantities for elliptic curves, and conjectured that the order of vanishing of the L-function $L(s)$ for the curve at $s \rightarrow 1$ equals to the rank.
This observation is the the now celebrated BSD Conjecture that bears their name; it is a Millennium Prize problem \cite{carlson2006millennium} and central to modern mathematics.

The above are but two of the countless examples where experimenting with pure data can lead to profound results.
They illustrate the importance of conjectures.
In theoretical research, finding the interesting problem is vital, and is time and again guided by the almost undefinable process of intuition.
G.~H.~Hardy's definition of mathematics is succinct: ``A mathematician, like a painter or a poet, is a maker of patterns.''.
Even a theoretical physicist or mathematical biologist, whose principle motivation comes from real world data and observations, would first distill the problem from Nature into a mathematical problem. Then it again becomes a mathematical game of pattern spotting, from graphs and plots, to formal symbols.
But here is an undeniable point: if there is one thing that AI can do better than humans, it is pattern recognition, especially when the data is in high dimension.

\subsection{Playing with Binary Sequences}
Let us perform the following simple experiment.
Given (i) the sequence \{ 0, 0, 1,  0, 0, 1,  0, 0, 1,  0, 0, 1,  0, 0, 1,  0, 0, 1 \} and asked what the next number is, any human would instantly say 0.
One way to describe it is the sequence of whether $n$ divides 3, for positive integers $n \in \mathbb{Z}_{>0}$.
Now, (ii) try \{ 0, 1, 1, 0, 1, 0, 1, 0, 0, 0, 1, 0, 1, 0, 0, 0, 1, 0, 1, 0, 0, 0, 1,
0, 0, 0, 0, 0, 1, 0, 1, 0, 0, 0, 0, 0, 1, 0, 0, 0, 1, 0, 1, 0, 0, 0, 
1, 0, 0, 0, 0, 0 \}.
An inspired person might, after some experimentation, conclude that the next number is 1; this is the sequence of PrimeQ, whether the $n$-th positive integer is prime or not.
Next, (iii) try \{ 1, 1, 1, 0, 1, 1, 1, 0, 0, 1, 1, 0, 1, 1, 1, 0, 1, 0, 1, 0, 1, 1, 1, 
0, 0, 1, 0, 0, 1, 1, 1, 0, 1, 1, 1, 0, 1, 1, 1, 0, 1, 1, 1, 0, 0, 1, 
1, 0, 0, 0, 1\}.
The person will be rather hard pressed to guess \footnote{I have tried this on the audience in numerous talks over the last few years.}.
This is the sequence of whether the $n$-th positive integer has a even (0) or odd (1) number of prime factors counted with multiplicity, a shifted version of the so-called M\"obius mu-function.
Uncovering its patterns would have incredible repercussions for mathematics: there are equivalent formulations of the Riemann Hypothesis in terms of this sequence.

What if we gave the sequence to AI?
For instance, what about a supervised machine learning algorithm?
In order to establish a reasonable training set, one could chose the following representation (and indeed the choice of representation is extremely important!).
Take one of the above infinite sequences $\{a_i\}_{i = 1,2,3, \ldots}$ and a sliding window of length $N$. In other words, consider a set of sequences $\{ 
\{a_i\}_{i=1,2,\ldots,N} \ , 
\{a_i\}_{i=2,3,\ldots,N+1} \ ,
\ldots
\{a_i\}_{i=k,k+1,\ldots,N+k-1}
\}$ for some $k$.
Here, $k$ will be taken to be sufficiently large (say 100,000) to create a decent data size, and $N$ will be taken to be sufficiently large (say 100) to give enough features.
After all, mathematical data is cheap!
We can then consider each of the finite sub-sequences as a single vector in $\mathbb{R}^N$ and label it by the next number outside the window:
\begin{equation}\label{window}
    \{
    (a_i)_{i=1,2,\ldots,N} \longrightarrow a_{N+1} \ , \quad
    (a_i)_{i=2,3,\ldots,N+1} \longrightarrow a_{N+2} \ , \quad
    \ldots \ , \quad
    (a_i)_{i=k,k+1,\ldots,N+k-1} \longrightarrow a_{N+k}
    \} \ .
\end{equation}
Note that we have chosen our sequences judiciously to standardize everything into a binary classification problem of binary vectors of dimension $N$.
The question is then: having seen $k$ labelled samples, how well will the ML algorithm predict on $k'$ unseen vectors.
This familiar supervised ML paradigm can be then compared to the human eye.

One can readily check that with the most basic ML algorithms suited for this problem, such as decision trees, support vector machines, relatively shallow feed-forward neural networks with ReLU activation functions, etc., on  Eq \eqref{window} applied to our three sequences.
For (i) one very quickly reaches 100\% accuracy for any of the ML strategies.
On (ii) one reaches \footnote{Note, in this case because of the increasing rarity of primes - approximately by a factor of $x/\ln(x)$ due to the Prime Number Theorem - we scale the window size accordingly.} about 80\%, while on (iii) one struggles to find any AI algorithm that would beat 50\%.
What this means is that (i), a trivial problem for the human eye, is as trivial for the AI; for (ii) the AI might be finding some version of the Sieve of Eratosthenes for checking PrimeQ; and for (iii) the AI is not beating a random guess. Of course, should one find an algorithm which {\it does}, then one might be well under way in finding a new approach toward the Riemann Hypothesis!

Many of the papers referenced in the introduction employ similar ideas to the previous paragraph, but to much more sophisticated situations.
Indeed, what (computable) mathematics, in the sense of a Turing Machine, does not fall into some version of Eq \eqref{window}?!
We can make the situation even more visual and suited to AI by wrapping each $N$ vector into an $m \times n = N$ matrix, which can be interpreted as a pixelated image where 1 is black and 0 is white.
For instance, suppose $N=100$, the first vector for situation (ii) can be wrapped into a $10 \times 10$ matrix, together with a label 1 (since 101 is prime): something like
$\begin{array}{c}\includegraphics[width=1cm]{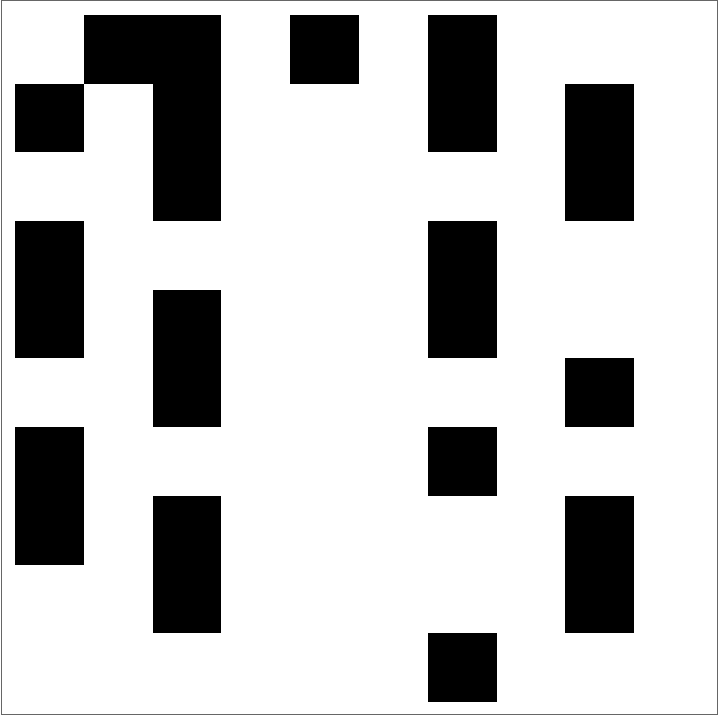}\end{array} \longrightarrow 1$.
Take, as a much more elevated example, the problem of computing a topological invariant for a manifold in algebraic geometry (which involves rather advanced calculations).
Yet, one can represent the manifold as a pixelated image by tensorizing the multi-degree information of the manifold as an algebraic variety \cite{He:2017aed,He:2018jtw}.
In a similar manner one could fashion any mathematical computation as an image recognition problem. Learning and gaining experience and intuition by a plethora of calculations - as mathematicians and theorists do during their careers - is in analogy to training a neural network. We could summarize this paradigm as
\begin{quote}
    {\it
    Bottom-up (and meta-) mathematics is language processing while top-down mathematics is image processing.
    }
\end{quote}

\subsection{The Birch Test}
The key steps to top-down mathematics are (1) identifying the problem and (2) identifying a strategy to attack the problem.
Both depend on experience, with a healthy dose of intuition.
While LLMs applied to the likes of Lean's MathLib \cite{de2015lean} are making baby-steps in (2), step (1) is formally known as ``Conjecture Formulation'', exemplified by the aforementioned cases of Gauss, Birch, Swinnerton-Dyer.
Can AI assist in telling a good conjecture from a useless one? Which patterns found from mathematical data lead to interesting as opposed to trivial mathematics?
More recently, this AI-guided ``Conjecture Formulation'' has been given much systematic thought \cite{wu2018can,He:2021oav,davies2021advancing,friederich2021scientific,raayoni2021generating,mishra2023mathematical,bauer2023mlfmf,davila2023advancements}.

In a six-month workshop in Cambridge in 2023 which I helped to co-organize \cite{ini2023}, together with Professor K.~Buzzard et al., we wanted to give some criteria on AI-driven theoretical discovery, and in particular on AI-assisted conjectures.
Since chatGTP has passed the Turing Test, we wanted to up the ante.
Inspired by a talk given by Birch \cite{birchtalk}, we called it the Birch Test \cite{he2024can}.
The AI-assisted discovery must be such that
\begin{description}
    \item[(A) Automaticity: ]
    it is completely made by AI from pattern-spotting, without any human intervention;
    \item[(I) Interpretability: ]
    any statements - conjectures or conclusions - be precise to a human mathematician, who cannot distinguish it from one given by a human colleague;
    \item[(N) Non-Triviality: ]
    it is non-trivial enough that the community of human experts will work on it.
\end{description}
To be fair, these are very stringent criteria (one needs to have a high bar in honour of Birch!) and so far no AI-assisted theoretical discovery has passed all three parts of the test. We now highlight with some examples where they succeed and fail.

Take the early experiments of obtaining topological invariants by deep-learning algebraic varieties \cite{He:2017aed,He:2018jtw}. They have been improved to $>99.9\%$ accuracy \cite{Erbin:2020tks}, which hint toward underlying and yet unknown structures in algebraic geometry that facilitate calculations without recourse to standard and computationally expensive sequence chasing.
These results suffer from the typical problem of deep neural networks: there is no interpretible formula one could extract. Thus, they fail Birch Test (I).
A better situation \cite{He:2019nzx} is where a support vector machine found a separation between simple and non-simple finite groups by plotting the Cayley multiplication tables. However, the hypersurface of separation is so complicated and furthermore deforms as more samples of groups were added that Birch Test (I) is still not passed.
The knot invariant relations found by saliency analyses \cite{davies2021advancing} and the Reidemeister moves untangled for extremely complicated knots \cite{Gukov:2020qaj}, though novel, interesting, and precise, were either readily proven or have not become sufficiently influential in the field; thus they fail Birch Test (N). Likewise, the continued fraction identities found by the Ramanujan machine \cite{raayoni2021generating}, or the physical conservation laws found by AI-Feynman \cite{udrescu2020ai} also belong to this category. Even the faster matrix multiplication algorithm found by DeepMind \cite{fawzi2022discovering} was shortly thereafter beaten by human researchers \cite{kauers2022fbhhrbnrssshk}.

The closest any AI-guided theoretical discovery to passing the Birch test so far, which has been precise enough and propelled a community of human experts to work on, is the Murmuration Conjectures in number theory \cite{He:2022pqn}. This discovery has passed (I) and for the first time passed (N). However, it fails (A) in that humans intervened in the process by digging under the hood: surprised by why AI was doing so well at distinguishing ranks of elliptic curves in the context of BSD, the researchers had to home in on a principal component analysis, and then look at the weight matrices in order to extract an unexpected formula.

\section{Prospectus}
The human theorist is not in danger of being put out of the job in the foreseeable future.
From the lack of a complete bottom-up MathLib database \cite{de2015lean} for all of mathematics, to the challenges LLMs would face given the vast search space of proof strategies even with such a database, to the exacting requirements of the Birch Test in top-down mathematics, we are far from automating theoretical discovery.

Nevertheless, it is undeniable that AI is beginning to play and will continue to play a pivotal role {\it in partnership} with the human mathematician and theoretical scientist.
In the nineteenth century, Gauss's intuitions alone were good enough to spot patterns that led to such profound results as the Prime Number Theorem.
In the twentieth century, computer experimentations were needed along-side the insights of Birch and Swinnerton-Dyer to raise the BSD Conjecture.
Now, in the twenty-first century, AI will work hand-in-hand with human experts to find new insights, conjectures, as well as strategies for derivations and proofs.


\section*{Acknowledgements}
I am most grateful to Dr.~Ananyo Bhattacharya,
Prof.~Alexander Kosyak and Dr.~Melissa Duncan for many valuable comments on the draft.
I would like to lend this opportunity to thank my many collaborators over the past seven years on AI-assisted mathematics, for the great fun and friendship: Daattavya Aggarwal, Laura Alessandretti, Guillermo 
Arias-Tamargo, Anthony Ashmore, Jiakang Bao, Andrea 
Baronchelli, Per Berglund, David Berman, Kieran Bull, Lucille 
Calmon, Hengyu Chen, Siqi Chen, Andrei Constantin, 
Pierre-Philippe Dechant, Rehan Deen, Stavros Garoufalidis, 
Elli Heyes, Edward Hirst, Johannes Hofscheier, Juan Ipiña, 
Vishnu Jejjala, Alexander Kasprzyk, Minhyong Kim, Shailesh 
Lal, Kyu-Hwan Lee, Seung-Joo Lee, Jianrong Li, Andre Lukas, 
Suvajit Majumder, Challenger Mishra, Gregg Musiker, Brent 
Nelson, Andrew Nestor, Thomas Oliver, Burt Ovrut, Toby 
Peterken, Stephen Pietromonaco, Andrey Pozdnyakov, Dmitrii 
Riabchenko, Diego Rodriguez-Gomez, Henrique Sá Earp, Max 
Sharnoff, Tomás Silva, Eldar Sultanow, Yuxuan Xiao, 
Shing-Tung Yau, and Zaid Zaz.

The resaerch is funded in part by STFC grant ST/J00037X/2 and the Leverhulme Trust for a project grant.

\section*{Competing interests}
The author declares no competing interests. 

\section*{Publisher’s note}
Springer Nature remains neutral with regard to jurisdictional claims in published maps and institutional affiliations.

\printbibliography

\end{document}